\documentclass[12pt]{amsart} 
\usepackage[mathscr]{eucal}
\usepackage{amsmath,amsfonts}

\parskip=\smallskipamount

\newtheorem{theorem}{Theorem}[section]

\newtheorem{lemma}[theorem]{Lemma}

\newtheorem{remark}[theorem]{Remark}
\newtheorem{problem}[theorem]{Problem}

\makeatletter
\@addtoreset{equation}{section}
\makeatother

\newcommand{\CC}{{\mathbb C}}

\newcommand{\FF}{{\mathbb F}}

\newcommand{\cB}{{\mathcal B}}

\newcommand{\cE}{{\mathcal E}}
\newcommand{\cF}{{\mathcal F}}
\newcommand{\cG}{{\mathcal G}}
\newcommand{\cH}{{\mathcal H}}

\newcommand{\cL}{{\mathcal L}}

\newcommand{\cR}{{\mathcal R}}
\newcommand{\cS}{{\mathcal S}}
\newcommand{\cT}{{\mathcal T}}
\newcommand{\cU}{{\mathcal U}}

\newdimen\expt
\def\boxit#1{\setbox0\hbox{$\displaystyle{#1}$}
      \hbox{\lower.4\expt
 \hbox{\lower3\expt\hbox{\lower\dp0
      \hbox{\vbox{\hrule height.4\expt
 \hbox{\vrule width.4\expt\hskip3\expt
      \vbox{\vskip3\expt\box0\vskip2\expt}%
 \hskip3\expt\vrule width.4\expt}\hrule height.4\expt}}}}}}
\begin{document}
\pagestyle{plain}

%\begin{flushright}
%  \it Date of this draft: \today
%\end{flushright}
\bigskip

\title %[] 
{A note on noncommutative interpolation} 
\author{T. Constantinescu} \author{J. L. Johnson} 

\address{Department of Mathematics \\
  University of Texas at Dallas \\
  Box 830688, Richardson, TX 75083-0688, U. S. A.}
\email{\tt tiberiu@utdallas.edu} 
\address{Department of Mathematics \\
  University of Texas at Dallas \\
  Box 830688, Richardson, TX 75083-0688, U. S. A   } 
\email{\tt jlj@utdallas.edu}

\begin{abstract}
In this paper we formulate and solve Nevanlinna-Pick and 
Carath\'eodory type problems for tensor algebras with data given
on the 
$N$-dimensional operator unit ball of a Hilbert space. 
We develop an approach based on the displacement structure
theory.
\end{abstract}

\maketitle

\section{Introduction}
Interpolation problems for bounded analytic functions were
studied quite intensively due to their many applications, for instance
to the wave propagation in layered media, circuit synthesis, or 
robust control. On the mathematical side, they are strongly related 
to dilation theory and selfadjoint extensions of symmetric operators.
Accordingly, it was some interest in generalizing the framework in which
similar problems could be formulated. In this paper we deal with 
interpolation problems for tensor algebras. 
Problems of this type were already considered in literature
(\cite{DP2}, \cite{Po2}). Several methods were developed
in these and other papers, 
each of them being of a specific interest. Here we consider similar
formulations for data on the $N$-dimensional unit ball
$$\cB _N(\cE )=\{Z=\left[
\begin{array}{ccc}
Z_1 & \ldots & Z_N
\end{array}
\right]\in \cL (\cE )^N \mid \sum _{k=1}^NZ^*_kZ_k<I_{\cE }\},$$
where $\cE $ is a Hilbert space. When considering this framework, we need
to introduce an evaluation of an element of the tensor algebra 
at a point of $\cB _N(\cE )$, as well as 
appropriate derivations on the algebra.
The goal of this paper is to define all these elements and formulate 
and solve Nevanlinna-Pick and Carath\'eodory type problems on 
$\cB _N(\cE )$. All of these are done by using a new approach based on 
the displacement structure theory, as suggested in \cite{CSK}. 
This approach turns out to be 
quite elementary and has other benefits, some of which are 
presented in our companion paper \cite{CJ1}.

The paper is organized as follows. In Section~2 we introduce a 
Szeg\"o type kernel for $\cB _N(\cE )$. Section~3 contains the 
description of our approach to interpolation
on $\cB _N(\cE )$ which is then illustrated by an application to a
Nevanlinna-Pick type problem. The last section introduces
some natural derivations on the tensor algebra and the 
formulation and solution of a Carath\'eodory type problem.
There is a suggestion here about the posibility of developing
a multidimensional calculus for tensor algebras that will
be pursued elsewhere. In this paper, which is part of 
\cite{Jo},  we restrict ourselves 
to the development of the displacement structure approach 
to interpolation for tensor algebras.

\section{Szeg\"o kernels}  

Let $\cE $ be a Hilbert space and let $\cL (\cE )$ denote 
the set of all bounded linear operators on $\cE $. If $N$ is a
positive integer and $Z=\left[
\begin{array}{ccc}
Z_1 & \ldots & Z_N
\end{array}
\right]$,
$W=\left[
\begin{array}{ccc}
W_1 & \ldots & W_N
\end{array}
\right]$
are two elements in $\cL (\cE )^N$, then we define 
$$(Z|W)=\sum _{k=1}^NZ^*_kW_k$$
and 
$$\cB _N(\cE )=\{Z=\left[
\begin{array}{ccc}
Z_1 & \ldots & Z_N
\end{array}
\right]\in \cL (\cE )^N \mid (Z|Z)<I_{\cE }\},$$
where $I_{\cE }$ denoted the identity operator on $\cE $.

We introduce a Szeg\"o type kernel on 
$\cB _N(\cE )$ by using some simple ideas from displacement
structure theory (and which, in the case $N=1$ and $\cE =\CC $, would give
the classical Szeg\"o kernel $K(z,w)=\frac{1}{1-\overline{z}w}$).
Let $Z_1$, $\ldots $, $Z_n$ be elements in $\cB _N(\cE )$ and consider
\begin{equation}\label{efuri}
F_k=\oplus _{l=1}^nZ^*_{l,k}, \quad \quad k=1,\ldots N,
\end{equation}
the diagonal matrix with the diagonal made of the $k$th components 
of $Z_1$, $\ldots $, $Z_n$. Also, define
\begin{equation}\label{u}
U=[\underbrace{\begin{array}{ccc}
I_{\cE } & \ldots & I_{\cE } 
\end{array}}_{n\,\,terms}]^*.
\end{equation}
It is easily seen that the so-called {\it displacement 
equation}
\begin{equation}\label{de}
A-\sum _{k=1}^NF_kAF^*_k=UU^*,
\end{equation}
admits a unique positive solution $A$. In fact, it is simple
to write the explicit form of the solution. Thus, let ${\FF }^+_N$
be the unital 
free semigroup with $N$ generators $1$, \ldots ,$N$. The empty
word is the identity element of ${\FF }^+_N$. The length 
of the word $\sigma $ is denoted by $|\sigma |$ and we consider
the lexicographic order on ${\FF }^+_N$.
We associate new Hilbert spaces to a Hilbert space $\cE $
by the following recursion: $\cE _0=\cE $ and for $k\geq 1$,
\begin{equation}\label{22}
\cE _k=\underbrace{\cE _{k-1}\oplus \ldots \cE _{k-1}}_{N\,\, terms}
=\cE ^{\oplus N}_{k-1}.
\end{equation}
Then $U_k=\left[F_{\sigma }U\right]_{|\sigma |=k}$
gives a bounded operator from $\cE _k$ into $\cE $, where 
$F_{\sigma }$ is a notation for the operator
$F_{i_1}\ldots F_{i_k}$ provided that $\sigma =i_1\ldots i_k$
(we set $F_{\emptyset }=I_{\cE }$).
One easily checks that 
$U^*_{\infty }=\left[U_k\right]_{k=0}^{\infty }$
is a bounded operator from $\oplus _{k=0}^{\infty }\cE _k$
into $\cE $ and the solution of \eqref{de} is given by the formula
$A=U^*_{\infty }U_{\infty }$. We introduce the notation
$Z^*_{\sigma }=Z^*_{i_1}\ldots Z^*_{i_k}$ for $\sigma =i_1\ldots i_k$;
$Z^*_{\sigma }$ should be distinguished from $(Z_{\sigma })^*$,
the adjoint of $Z_{\sigma }$. Also define $L(Z)=
\left[Z_{\sigma }\right]_{|\sigma |=0}^{\infty }$
for $Z$ in $\cB _N(\cE )$. Then the solution of \eqref{de}
can be written in the form
$$A=\left[L(Z_j)L(Z_k)^*\right]_{j,k=1}^{n}.$$
This formula suggests to introduce the positive definite kernel
\begin{equation}\label{szego}
K(Z,W)=L(Z)L(W)^*,\quad \quad Z,W \in \cB _N(\cE ),
\end{equation}
as another generalization of the classical Szeg\"o kernel. 
Also, when
$\cE =\CC$ and $N>1$, we obtain  that 
\begin{equation}\label{fara}
K(Z,W)=(1-(Z|W))^{-1},
\end{equation}
which is a positive definite kernel on the unit ball in the complex 
$N$-dimensional space $\CC ^N$ that was studied quite intensively 
in recent years (\cite{Po2}, \cite{Ar}). Most notably, the 
kernel \eqref{fara} has a universality property
with respect to the Nevanlinna-Pick problem, as explained in
\cite{AM}. Note that for $\dim \cE >1$, the kernel 
\eqref{fara}
is no longer positive definite.

\section{Nevanlinna-Pick interpolation}

In this section 
we formulate and solve a Nevanlinna-Pick interpolation
problem for the noncommutative analytic Toeplitz algebras
introduced in \cite{Po1}. These algebras represent a multidimensional
generalization of the classical Toeplitz algebra assocaited to 
the Hardy space $H^{\infty }$. One reason for their study is that
they are Banach algebras containing the tensor algebra.
The associative tensor algebra $\cT (\cH)$ generated by the 
complex vector space
$\cH=\CC ^N$ is defined by the algebraic direct sum
$$\cT (\cH)=\oplus _{k\geq 0}\cH^{\otimes k},$$
where 
$\cH^{\otimes k}=
\underbrace{\cH\otimes \ldots \otimes \cH}_{k\,\, factors}$
is the $k$-fold algebraic tensor product of $\cH$ with itself.
Also, the Hilbert direct sum
$$\cF (\cH)=\oplus _{k\geq 0}\cH^{\otimes k}$$ 
is the full Fock space associated to $\cH $
(on each $\cH^{\otimes k}$, $k>1$, 
we   
consider the tensor Hilbert space structure induced by the 
Euclidean norm on $\cH={\CC}^N$; also $\cH^{\otimes k}=\CC$ and 
$\cH^{\otimes 1}\cH$ - see \cite{Pa} for more details on Fock space
constructions).
The noncommutative Toeplitz algebra, \cite{Po1}, 
can be identified with the set of those $\phi \in \cF (\cH)$
such that 
$$\sup \{\|\phi \otimes p\|_{\cF (\cH)}\mid p\in \cT (\cH),
\|p\|_{\cF (\cH)}\leq 1\}<\infty .$$

We notice that each $\cH^{\otimes k}$ can be identified 
with the Hilbert space $\cH _k$ given by \eqref{22}
and the noncommutative Toeplitz algebra
is isometrically isomorphic to the algebra $\cU _{\cT}(\cH)$
of upper triangular operators 
$T=[T_{ij}]_{i,j=0}^{\infty }\in \cL (\oplus _{k=0}^{\infty }\cH_k)$
with the property that for $i\leq j$ and $i,j\geq 1$, 
\begin{equation}\label{condition}
T_{ij}=T_{i-1,j-1}^{\oplus N},
\end{equation}
where for an operator $T$, $T^{\oplus N}=\underbrace{T\oplus
\ldots \oplus T}_{N\,\, terms}$. We also use the notation 
$\mbox{diag}[T]$ for the direct sum of a certain number of copies of $T$.
Denote by $\cS (\cH)$ the Schur class of all contractions 
in $\cU _{\cT}(\cH)$. 
Given a Hilbert space $\cE $, we can introduce the algebra
$\cU _{\cT}(\cH ,\cE)$ to be the set of all 
$T=[T_{ij}]_{i,j=0}^{\infty }\in \cL (\oplus _{k=0}^{\infty }\cE_k)$
satisfying \eqref{condition}. The corresponding Schur class of all 
contractions in $\cU _{\cT}(\cH ,\cE)$ is denoted by $\cS (\cH ,\cE)$.
Note that $\cE _k$ can be identified with $\cH _k\otimes \cE$, 
which justifies our notation.

Since the noncommutative Toeplitz algebras were viewed as generalizations
of the classical Toeplitz algebra, it was quite natural to study
bounded interpolation problems in this setting. To that end, 
a "point evaluation" was introduced and studied in \cite{Po1}, \cite{Ar}
and references therein. At about the same time, bounded
interpolation problems were studied for the algebra of upper
triangular operators in \cite{BGK}, \cite{DD}, \cite{SCK}
(see \cite{Con} for details and other related references)
and a point evaluation was introduced in this setting too.  
Since $\cU _{\cT}(\cH ,\cE)$ is an algebra of upper triangular
operators we can use the later approach as follows:
for $T\in \cU _{\cT}(\cH ,\cE)$ and $Z\in \cB _N(\cE )$
define the operator 
\begin{equation}\label{evaluare}
T(Z)=P_{\cE }TL(Z)^*,
\end{equation}
where $P_{\cE }$ denotes the orthogonal projection 
of $\oplus _{k=0}^{\infty }\cE_k$ onto $\cE (=\cE _0)$.
The basic property that qualifies this operator as a point evaluation 
is given by the following result.

\begin{lemma}\label{L31}
If $T\in \cU _{\cT}(\cH ,\cE)$, then 
$TL(Z)^*=\mbox{diag}[T(Z)]L(Z)^*.$
\end{lemma}

\begin{proof}
We write $L(Z)=\left[\begin{array}{cccc}
I & \tilde Z^*_1 & \tilde Z^*_2 & \ldots 
\end{array}\right]$, 
where $\tilde Z^*_k=\left[\begin{array}{c}
Z^*_{\sigma }\end{array}\right]_{|\sigma |=k}$.
Due to the properties of the lexicographic order
we deduce that $\tilde Z^*_k=
\left[\begin{array}{cc}
Z^*_1\tilde Z^*_{k-1} & \ldots Z^*_N\tilde Z^*_{k-1}
\end{array}\right]$. Therefore the $k$th block entry of $TL(Z)^*$
is
$$\begin{array}{rcl}
\sum _{l=0}^{\infty }T_{k,k+l}\tilde Z^*_{k+l} &=&
\sum _{l=0}^{\infty }T^{\oplus N}_{k-1,k+l-1}
\left[\begin{array}{c}
\tilde Z_{k-1}Z_1 \\
\vdots \\
\tilde Z_{k-1}Z_N
\end{array}\right] \\
 & & \\
 &=&
\left[\begin{array}{c}
\sum _{l=0}^{\infty }T_{k-1,k+l-1}\tilde Z_{k-1}Z_1 \\
\vdots \\
\sum _{l=0}^{\infty }T_{k-1,k+l-1}\tilde Z_{k-1}Z_N
\end{array}\right].
\end{array}
$$
This formula can be used in an inductive argument in order 
to conclude the proof.
\end{proof}

We can now formulate the following Nevanlinna-Pick type
problem:

\begin{problem}\label{NP}
Determine for which $Z_1$, $\ldots $, $Z_n$ in $\cB _N(\cE )$ and 
$B_1$, $\ldots $, $B_n$ in $\cL (\cE)$
there is a $T\in \cS (\cH ,\cE)$ such that
$T(Z_k)=B_k$, $k=1,\ldots n.$
\end{problem}

This problem can be solved by using the methods in \cite{Po2}, but 
here we indicate a new, elementary approach based on the displacement
structure theory. More precisely, we use a result announced in 
\cite{CSK} and proved in details in \cite{CJ1}, that gives the 
solution of a so-called scattering experiment associated to the 
data of the Nevanlinna-Pick problem. 
We will show that this data can be encoded by a
displacement equation of the following type:

\begin{equation}\label{newdisplacement}
A-\sum _{k=1}^NF_kAF_k^*=GJG^*,
\end{equation}
where $F_k\in \cL(\cG)$, $k=1,\ldots ,N$, are given contractions
on the Hilbert
space $\cG$. Also $G=\left[\begin{array}{cc}
U & V
\end{array}\right]\in \cL (\cE ^2,\cG)$ and $J=\left[
\begin{array}{cc}
I_{\cE } & 0 \\
0 & -I_{\cE } 
\end{array}\right]$.
The wave operators associated to \eqref{newdisplacement} 
are introduced by the formulae:
$U^*_{\infty }=\left[U_k\right]_{k=0}^{\infty }$,
$V^*_{\infty }=\left[V_k\right]_{k=0}^{\infty }$, where
$U_k=\left[F_{\sigma }U\right]_{|\sigma |=k}:\cE _k\rightarrow \cG$,  
$V_k=\left[F_{\sigma }V\right]_{|\sigma |=k}:\cE _k\rightarrow \cG$,
and $\cH=\CC ^N$. 
We will assume that both $U_{\infty } $ and $V_{\infty } $ are bounded
and also that $\lim _{k\rightarrow \infty }
\sum _{|\sigma |=k}\|F^*_{\sigma }g\|=0$
for all $g\in \cG $. Under these assumptions we
deduce that
\eqref{newdisplacement} has a unique solution given by
\begin{equation}\label{solution}
A=U^*_{\infty }U_{\infty }-V^*_{\infty }V_{\infty }. 
\end{equation}

\begin{theorem}\label{T32}
The solution \eqref{solution}
of the displacement equation \eqref{newdisplacement} 
is positive if and only if there exists 
$T\in \cS (\cH ,\cE)$ such that $V_{\infty }=TU_{\infty }$.
\end{theorem}
 
\begin{proof}
For sake of completeness
we indicate the main ideas of the proof. This will also show that 
our approach is quite elementary.

Assume $A=U^*_{\infty }U_{\infty }-V^*_{\infty }V_{\infty }\geq 0$ 
and let $A=LL^*$ be 
a factorization of $A$ with $L\in \cL (\cF ,\cG )$ for some 
Hilbert space $\cF $. From \eqref{newdisplacement} we deduce
that 
$$LL^*+VV^*=\sum _{k=1}^NF_kLL^*F^*_k+UU^*.$$
In matrix form, 
\begin{equation}\label{base}
\left[\begin{array}{cc}
L & V
\end{array}\right]
\left[\begin{array}{c}
L^* \\
V^*
\end{array}\right]=
\left[\begin{array}{cccc}
F_1L & \ldots & F_NL & U
\end{array}\right]
\left[\begin{array}{c}
L^*F^*_1 \\
\vdots \\
L^*F^*_N \\
 U^*
\end{array}\right].
\end{equation}  
Define $A^*=\left[\begin{array}{cc}
L & V
\end{array}\right]$ and 
$B^*=\left[\begin{array}{cccc}
F_1L & \ldots & F_NL & U
\end{array}\right]$, then we deduce from \eqref{base}
that there exists an unitary operator
$\theta _0\in \cL (\overline{\cR (B)},\overline{\cR (A)})$ such that
$A=\theta _0B.$
It follows that there exist Hilbert spaces 
$\cR _1$, $\cR _2$, and an unitary extension 
$\theta \in \cL (\cF ^{\oplus N}\oplus \cE \oplus \cR _1,
\cF \oplus \cE \oplus \cR _2)$ of $\theta _0$, hence this extension  
satisfies the relation 
\begin{equation}\label{relation}
\left[\begin{array}{c}
A \\
0_{\cR _2}
\end{array}\right]
=\theta \left[\begin{array}{c}
B \\
0_{\cR _1}
\end{array}\right].
\end{equation}

\noindent
Let $\theta _{ij}$, $i\in \{1,2,3\}$, $j\in \{1,2,\ldots ,N+2\}$,
be the matrix coefficients of $\theta $. It is convenient
to rename some of these coefficients. Thus, we set
$$\begin{array}{cccc}
X_k=\theta _{1k}, & k=1,\ldots N, & \quad & Z=\theta _{1,N+1}, \\
 & &  & \\
Y_k=\theta _{2k}, & k=1,\ldots N, & \quad & W=\theta _{2,N+1}.
\end{array}$$
From \eqref{relation} we deduce that
$$L^*=\sum _{k=1}^NX_kL^*F^*_k+ZU^*$$
and 
$$V^*=\sum _{k=1}^NY_kL^*F^*_k+WV^*.$$
By induction we deduce that 
\begin{equation}\label{key}
V^*=WU^*+\sum _{k=1}^N\sum _{|\sigma |=0}^nY_kX_{\sigma }ZU^*F_{k\sigma }^*
+\sum _{|\tau |=n+1}Q_{\tau }L^*F_{\tau }^*,
\end{equation}
where $Q_{\tau }$ are monomials of length $|\tau |$ in the variables
$X_1$, $\ldots $,$X_N$, $Y_1$, $\ldots $, $Y_N$.
Since $\theta $ is unitary it follows that all $Q_{\tau }$ are contractions.

We define $T_{00}=W$ and for $j>0$,
$$T_{0j}=[Y_kX_{\sigma }Z]_{|\sigma |=j-1;k=1,\ldots ,N}.$$
Then we define $T_{ij}$, $i>0$, $j\geq i$, by the 
formula \eqref{condition} and $T_{ij}=0$ for $i>j$. It can be checked 
that $T=[T_{ij}]_{i,j=0}^{\infty }$ belongs to $\cS (\cH ,\cE)$. 
Also, since 
$\lim _{k\rightarrow \infty }
\sum _{|\sigma |=k}\|F^*_{\sigma }g\|=0$
for all $g\in \cG $, we deduce from \eqref{key}
that $V_{\infty }=TU_{\infty }$.
\end{proof}

We can now give a solution to Problem~\eqref{NP}.

\begin{theorem}\label{T33}
Let $Z_1$, $\ldots $, $Z_n$ be distinct elements in $\cB _N(\cE )$ and 
$B_1$, $\ldots $, $B_n$ in $\cL (\cE)$. Then 
there is a $T\in \cS (\cH ,\cE)$ such that
$T(Z_k)=B_k$, $k=1,\ldots n$, if and only if the Pick matrix
$$P=\left[\begin{array}{c}
L(Z_j)\mbox{diag}\left[I-B_j^*B_k\right]L(Z_k)^*
\end{array}\right]_{j,k=1}^n$$
is positive.
\end{theorem}

\begin{proof}
First assume that $P$ is positive. Define the operators $F_k$, 
$k=1,\ldots ,N$, by the formula \eqref{efuri}. Let
$V=\left[\begin{array}{ccc}
B_1 & \ldots & B_n 
\end{array}\right]^*$ and set $G=\left[\begin{array}{cc}
U & V 
\end{array}\right]$, where $U$ was defined by \eqref{u}.
Then the unique solution of the displacement equation
$$
A-\sum _{k=1}^NF_kAF_k^*=GJG^*
$$
is $P$, which is positive. By Theorem~\ref{T32}, there is
$T\in \cS (\cH ,\cE )$ such that $V_{\infty }=TU_{\infty }$.
We can check that
$$U_{\infty }=\left[\begin{array}{ccc}
L(Z_1)^* & \ldots & L(Z_n)^* 
\end{array}\right],$$
$$V_{\infty }=\left[\begin{array}{ccc}
\mbox{diag}\left[B_1\right]L(Z_1)^* & \ldots & 
\mbox{diag}\left[B_n\right]L(Z_n)^* 
\end{array}\right].$$
From $V_{\infty }=TU_{\infty }$ we deduce that 
$TL(Z_k)^*=\mbox{diag}\left[B_k\right]L(Z_k)^*$, $k=1,\ldots ,n$.
By Lemma~\ref{L31} we deduce that $\mbox{diag}\left[B_k\right]L(Z_k)^*=
\mbox{diag}\left[T(Z_k)\right]L(Z_k)^*$, which implies that 
$T(Z_k)=B_k$, $k=1, \ldots n$.

Conversely, assume that there is 
a $T\in \cS (\cH ,\cE)$ such that
$T(Z_k)=B_k$, $k=1,\ldots n$. Then , with the previous notation, 
we deduce that $V_{\infty }=TU_{\infty }$ and so, 
$$P=U^*_{\infty }U_{\infty }-V^*_{\infty }V_{\infty }
=U^*_{\infty }(I-T^*T)U_{\infty }\geq 0$$
since $T$ is a contraction.
\end{proof}

If $Z_k=\mbox{diag}\left[z_k\right]$ for some $z_k\in \CC$, 
$k=1, \ldots n$, then
Problem~\ref{NP} reduces to the Nevanlinna-Pick type
formulated in \cite{DP2} and \cite{Po2}. 
In this case, Theorem~\ref{T33} reduces to the results   
given in \cite{DP2} and \cite{Po2} for the solution of this problem. 

\section{Differentiation on $\cU _{\cT}(\cH ,\cE)$}

We introduce several derivations on $\cU _{\cT}(\cH ,\cE)$.
This might be of interest in intself and also it allows us to
formulate interpolation problems involving higher derivates.
The motivation for the definitions comes from similar constructions
for the algebra of upper triangular operators given in 
\cite{SCK} (see \cite{Con} for more details and additional references).

Let $Z\in \cB _N(\cE )$, $Z=\left[
\begin{array}{ccc}
Z_1 & \ldots & Z_N
\end{array}
\right]$. We introduce the lower triangular operators
$F_k^{(l)}=\left[X_{ij}\right]\in \cL (\oplus _{k=0}^l\cE _k)$
with $X_{00}=Z_k^*$, $X_{10}=E_k$, $X_{i0}=0$ for $i>1$, and 
otherwise $X_{ij}=X^{\oplus N}_{i-1,j-1}$.
Also define

\begin{equation}\label{uri}
\begin{array}{cccc}
U & = & [I_{\cE } 
& \underbrace{ 0 \,\,\,\, \ldots \,\,\,\, 0 }_{N+\ldots +N^l \,\, terms}]^*.
\end{array}
\end{equation}

Then $U^*_{\infty }=\left[F_{\sigma }^{(l)}U
\right]_{|\sigma |=0}^{\infty }$ and notice that
$$U_{\infty }=\left[\begin{array}{cc}
L(Z)^* & \left[L(\sigma ,Z)\right]^*\end{array}\right]_{|\sigma |=1}^l,$$
where $L(\sigma ,Z)$ are well-defined operators whose form follows
from the previous relation. 
We define the partial derivatives
\begin{equation}\label{partiale}
D_{\sigma }T_Z=P_{\cE }TL(\sigma ,Z)^*, \sigma \in 
{\FF }^+_N;
\end{equation}
also set $D_{\emptyset }T_Z=T(Z)$.

\begin{remark}
{\rm
We notice that 
$$F_{\sigma }^{(1)}=
\left[\begin{array}{cc}
Z^*_{\sigma } & 0 \cr
P_{\sigma } & (Z^*_{\sigma })^{\oplus N} 
\end{array}\right],
$$
where $P_{\sigma } $ are defined by the recursions:
$P_k=E_k$, $k=1,\ldots ,N$, and 
\begin{equation}\label{rec}
P_{k\tau }=E_kZ^*_{\tau }+(Z^*_k)^{\oplus N}P_{\tau },
\quad k=1,\ldots ,N, \,\,\mbox{and}\,\, |\tau |>0.
\end{equation}
We deduce from here that 
$L_k(Z)=L(Z)A_k$ for some densely defined unbounded
operator $A_k$. We can formally write
$$D_{k}T_Z=P_{\cE }TA_k^*L(Z)^*,$$
so that, at least formally, we can define $P_kT=TA^*$. 
Of course, this operator cannot be bounded, but its 
first row is a well defined bounded operator,
$S^{(1)}=\left[S_{0k}\right]_{k=0}^{\infty }$.
Therefore, we can introduce a densely defined operator $D_kT$
on $\oplus _{k=0}^{\infty }\cE _k$, such that 
$(D_kT)_{0k}=S_{0k}$, $k\geq 0$, and otherwise
$(D_kT)_{ij}=(D_kT)_{i-1,j-1}$.
Then $\delta _k(T)=D_kT$ gives, formally, a derivation 
on $\cU _{\cT}(\cH ,\cE)$. More details about these derivations
are considered in \cite{CJ2}.
}
\end{remark}

A Carath\'eodory type problem can be formulated
as follows.

\begin{problem}\label{C}
Given $Z\in \cB _N(\cE )$ and $l$ a positive integer, determine for which 
$B_k\in \cL (\cE _k,\cE)$, $0\leq k\leq l$,
there is a $T\in \cS (\cH ,\cE)$ such that
$\left[D_{\sigma }T_Z\right]_{|\sigma |=k}=B_k$ for $0=1,\ldots l.$
\end{problem}

The solution of Problem~\ref{C} can be obtained 
by a construction
similar to the one involved in the proof of Theorem~\ref{T33}.
Thus, let $V=\left[\begin{array}{ccc}
B_1 & \ldots & B_l
\end{array}\right]^*$
and set $G=\left[\begin{array}{cc}
U & V \end{array}\right]$, where $U$ is defined by \eqref{uri}.
Then the unique solution of the displacement equation
$$A-\sum _{k=1}^NF_k^{(l)}A(F_k^{(l)})^*=GJG^*$$
is $A=U_{\infty }^*U_{\infty }-V_{\infty }^*V_{\infty }$, where 
we use the notation involved in the statement of Theorem~\ref{T32}.
We obtain the following result.

\begin{theorem}\label{T43}
Let $Z\in \cB _N(\cE )$ and 
$B_k\in \cL (\cE _k,\cE)$, $0\leq k\leq l$ be given. Then 
there is a $T\in \cS (\cH ,\cE)$ such that
$\left[D_{\sigma }T_Z\right]_{|\sigma |=k}=B_k$, $k=0,\ldots l,$
if and only if the matrix
$U_{\infty }^*U_{\infty }-V_{\infty }^*V_{\infty }$
is positive.
\end{theorem}

\begin{proof}
This is similar to the proof of Theorem~3.4 and we can omit the details.
\end{proof}

\noindent
We notice that for $Z=0$ we deduce the Corollary~5.2 in \cite{Po2}.
A more explicit form for $U_{\infty }^*U_{\infty }-V_{\infty }^*V_{\infty }$
can be obtained but that is only a notational matter.

We conclude this section by introducing another possible formulation
of a Carath\'eodory type problem in this setting. Thus, we can introduce
a total derivate of order $k$ at $Z$ by the formula:
\begin{equation}\label{totala}
D^kT_Z=\sum _{|\sigma |=k}D_{\sigma }T_Z, \quad k=1, \ldots ,
\end{equation}
and by convention $D^0T_Z=T(Z)$. A corresponding Carath\'eodory type
problem can be formulated as follows.

\begin{problem}\label{C2}
Given $Z\in \cB _N(\cE )$ and $l$ a positive integer, determine for which 
$B_k\in \cL (\cE)$, $0\leq k\leq l$,
there is a $T\in \cS (\cH ,\cE)$ such that
$D^kT_Z=B_k$ for $k=0,\ldots l.$
\end{problem}
 
We can show that this problem can be also solved by using the 
displacement structure approach. Thus, let $l$ be a positive integer
and $Z\in \cB _N(\cE )$, $Z=\left[
\begin{array}{ccc}
Z_1 & \ldots & Z_N
\end{array}
\right]$. Then, for $1\leq k\leq N$,  
$$TF_k^{(l)}=\left[\begin{array}{ccccc}
Z^*_k & 0 & \ldots & & \\
I_{\cE} & Z^*_k & 0 & \ldots & \\
0 & I_{\cE} & Z^*_k  & \ldots \\
\vdots & & & \ddots &   \\
0 & \ldots & 0 & I_{\cE } & Z^*_k
\end{array}\right]$$
and

\begin{equation}\label{uri2}
\begin{array}{cccc}
U & = & [I_{\cE } 
& \underbrace{ 0 \,\,\,\, \ldots \,\,\,\, 0 }_{l \,\,terms}]^*.
\end{array}
\end{equation}
The associated $U_{\infty }^*=
\left[TF^{(l)}_{\sigma }U\right]_{|\sigma |=0}^{\infty }$
can be written in the form 
$$U_{\infty }=\left[
\begin{array}{cccc}
L(Z) & M_1(Z) &\ldots & M_l(Z)
\end{array}\right],
$$
for some operators $M_k(Z)$, $k=1$, $\ldots $, $l$.

\begin{lemma}\label{L43}
$$D^kT_Z=P_{\cE }TM_k(Z).$$
\end{lemma}

\begin{proof}
We prove by induction that
$$F_{\sigma }^{(l)}=
\left[\begin{array}{ccccc}
Z^*_{\sigma } & 0 & \ldots & & \\
P_{1}(\sigma ) & (Z^*_{\sigma })^{\oplus N} & 0 & \ldots & \\
P_{2}(\sigma ) & (P_{1,\sigma })^{\oplus N} & 
(Z^*_{\sigma })^{\oplus N^2} & \ldots \\
\vdots & & & \ddots &   \\
P_{l}(\sigma ) & & \ldots   & & (Z^*_{\sigma })^{\oplus N^l}
\end{array}\right]$$
and with $P_{0}(\sigma )=Z^*_{\sigma },$
\begin{equation}\label{relu}
P_{j}(k\sigma )=E_k^{\oplus N^{j-1}}P_{j-1}(\sigma )
+(Z^*_{\sigma })^{\oplus N^j}P_{j}(\sigma ),\quad \quad k=1, \ldots ,N.
\end{equation}
Similarly,
$$TF_{\sigma }^{(l)}=
\left[\begin{array}{ccccc}
Z^*_{\sigma } & 0 & \ldots & & \\
Q_{1}(\sigma ) & Z^*_{\sigma } & 0 & \ldots & \\
Q_{2}(\sigma ) & Q_{1}(\sigma ) & 
Z^*_{\sigma } & \ldots \\
\vdots & & & \ddots &   \\
Q_{l}(\sigma ) & & \ldots   &  Q_{1}(\sigma ) & Z^*_{\sigma }
\end{array}\right]$$
and with $Q_0(\sigma )=Z^*_{\sigma },$
\begin{equation}\label{reld}
Q_{j}(k\sigma )=Q_{j-1}(\sigma )
+Z^*_{\sigma }Q_{j}(\sigma ),\quad \quad k=1,\ldots ,N.
\end{equation}
Now, each $P_{j}(\sigma )$ is a column matrix with 
$N^j$ entries $P^s_{j}(\sigma )$.
From \eqref{relu} and \eqref{reld} it follows that 
$$\sum _{s=1}^{N^j}P^s_j(\sigma )=Q_j(\sigma )$$
and this implies the required formula.
\end{proof}

This remark allows us tho solve Problem~\ref{C2} in the same way
we solved Problem~\ref{C}. Thus, define
$V=\left[\begin{array}{ccc}
B_1 & \ldots & B_l
\end{array}\right]^*$
and set $G=\left[\begin{array}{cc}
U & V \end{array}\right]$, where $U$ is defined by \eqref{uri2}.
Then the unique solution of the displacement equation
$$A-\sum _{k=1}^NF_k^{(l)}A(F_k^{(l)})^*=GJG^*$$
is $A=U_{\infty }^*U_{\infty }-V_{\infty }^*V_{\infty }$, where again
we use the notation involved in the statement of Theorem~\ref{T32}.
With the same proof as that of Theorem~\ref{T43} 
we obtain the following result.

\begin{theorem}\label{T44}
Let $Z\in \cB _N(\cE )$ and 
$B_k\in \cL (\cE)$, $0\leq k\leq l$ be given. Then 
there is a $T\in \cS (\cH ,\cE)$ such that
$D^kT_Z=B_k$ for $k=0,\ldots l$
if and only if the matrix
$U_{\infty }^*U_{\infty }-V_{\infty }^*V_{\infty }$
is positive.
\end{theorem}


\begin{thebibliography}{99}
\frenchspacing

\bibitem{AM}
J.~Agler and J.~E.~McCarthy,
{\em Complete Nevanlinna-Pick kernels}, 
preprint 1997.

\bibitem{Ar}
W.~Arveson, 
Subalgebras of $C^*$-algebras III:
Multivariable operator theory,
{\em Acta Math.}, {\bf 181}(1998), 476--514. 

%\bibitem{BTV}
%J.~A.~Ball, T.~T.~Trent, and V.~Vinnikov,
%Interpolation and commutant lifting for multipliers
%on reproducing kernel Hilbert spaces, preprint 2000.

\bibitem{BGK}
J.~A.~Ball, I.~Gohberg, and M.~A.~Kaashoek,
Nevanlinna-Pick interpolation problem for time-varying
input-output maps: The discrete case, in 
{\em Operator Theory: Advances and Applications},
Vol. 56, Birkh\"auser, 1992, pp. 1--51.



\bibitem{Con}
T.~Constantinescu,
{\em Schur Parameters, Factorization and Dilation Problems},
Birkh\"auser, 1996.

\bibitem{CJ1}
T.~Constantinescu, and J.~L.~Johnson,
Tensor algebras, displacement structure, and the Schur algorithm, 
preprint 2001.


\bibitem{CJ2}
T.~Constantinescu, and J.~L.~Johnson,
Calculus on the tensor algebra, in preparation.


\bibitem{CSK}
T.~Constantinescu, A.~H.~Sayed, and T.~Kailath,
Inverse scattering experiments, structured matrix 
inequalities, and tensor algebras, {\em Linear Alg. Appl.}, to appear.

\bibitem{DP1}
K.~R.~Davidson and D.~R.~Pitts,
The algebraic structure of non-commutative analytic Toeplitz
algebras, {\em Math. Ann.}, {\bf 311}(1998), 275--303.

\bibitem{DP2}
K.~R.~Davidson and D.~R.~Pitts,
Nevanlinna-Pick interpolation for noncommutative 
analytic Toeplitz algebras,

\bibitem{DD}
P.~Dewilde and H.~Dym, 
Interpolation for upper triangular operators, in
{\em Operator Theory: Advances and Applications},
Vol. 56, Birkh\"auser, 1992, pp. 153--260.


\bibitem{Jo}
J.~L.~Johnson, 
{\em Tensor Algebras, Displacement Structure, and some Classes
of Stochastic Processes}, Dissertation, University of Texas at Dallas, 
in preparation.

\bibitem{KS}
T.~Kailath and A.~H.~Sayed,
Displacement structure: theory and applications,
{\em SIAM Rev.}, {\bf 37}(1995), 297--386. 

\bibitem{Pa}
K.~R.~Parthasarathy,
{\em An Introduction to Quantum Stochastic Calculus},
Birkh\"auser, 1992.

\bibitem{Po1}
G.~Popescu,
Multi-analytic operators on Fock spaces,
{\em Math. Ann.}, {\bf 303}(1995), 31--46.

\bibitem{Po2}
G.~Popescu,
Interpolation problems in several variables,


\bibitem{SCK}
A.~H.~Sayed, T.~Constantinescu, and T.~Kailath,
Lattice structures for time variant interpolation problems,
in {\em Proc. IEEE Conf. Decision and Contr.}, 
vol. 1, Tucson, 1992, pp. 116--121.


\end{thebibliography}
\end{document}